\documentclass[12pt]{article}
\usepackage[dvips]{epsfig}
\usepackage{amsthm,amsmath,amssymb}
\usepackage{latexsym}
\usepackage{placeins}
\usepackage{color}
\usepackage{hyperref}
\hypersetup{
     colorlinks=true,
     linkcolor=blue,
     filecolor=blue,
     citecolor = black,      
     urlcolor=cyan,
     }
\definecolor{blue}{rgb}{0,0,0.9}
\definecolor{red}{rgb}{0.9,0,0}
\definecolor{green}{rgb}{0,0.50,0.10}
\definecolor{violet}{rgb}{0.5804,0.0000,0.8275}

\def\@themcountersep{}

\newtheorem{THEO}{Theorem}[section]
\newtheorem{ALGo}[THEO]{Algorithm}
\newtheorem{CONJ}[THEO]{Conjecture}
\newtheorem{COND}[THEO]{Condition}
\newtheorem{ASSUMP}[THEO]{Assumption}
\newtheorem{CORO}[THEO]{Corollary}
\newtheorem{DEFI}[THEO]{Definition}
\newtheorem{EXAMP}[THEO]{Example}
\newtheorem{FACT}[THEO]{Fact}
\newtheorem{HYPO}[THEO]{Hypothesis}
\newtheorem{LEMM}[THEO]{Lemma}
\newtheorem{PROB}[THEO]{Problem}
\newtheorem{PROP}[THEO]{Proposition}
\newtheorem{REMA}[THEO]{Remark}
\newcommand{\theo}{\begin{THEO}}
\newcommand{\algo}{\begin{ALGo} \rm}
\newcommand{\cond}{\begin{COND} \rm}
\newcommand{\assump}{\begin{ASSUMP} \rm}
\newcommand{\conj}{\begin{CONJ}}
\newcommand{\coro}{\begin{CORO}}
\newcommand{\defi}{\begin{DEFI} \rm}
\newcommand{\examp}{\begin{EXAMP} \rm}
\newcommand{\fact}{\begin{FACT}}
\newcommand{\hypo}{\begin{HYPO} \rm}
\newcommand{\lemm}{\begin{LEMM}}
\newcommand{\prob}{\begin{PROB} \rm}
\newcommand{\prop}{\begin{PROP}}
\newcommand{\rema}{\begin{REMA} \rm}
\newcommand{\etheo}{\end{THEO}}
\newcommand{\ealgo}{\end{ALGo}}
\newcommand{\econd}{\end{COND}}
\newcommand{\eassump}{\end{ASSUMP}}
\newcommand{\econj}{\end{CONJ}}
\newcommand{\ecoro}{\end{CORO}}
\newcommand{\edefi}{\end{DEFI}}
\newcommand{\eexamp}{\end{EXAMP}}
\newcommand{\efact}{\end{FACT}}
\newcommand{\ehypo}{\end{HYPO}}
\newcommand{\elemm}{\end{LEMM}}
\newcommand{\eprob}{\end{PROB}}
\newcommand{\eprop}{\end{PROP}}
\newcommand{\erema}{\end{REMA}}

\def\0{\mbox{\bf 0}}
\def\1{\mbox{\bf 1}}
\def\2{\mbox{\bf 2}}
\def\3{\mbox{\bf 3}}
\def\4{\mbox{\bf 4}}
\def\5{\mbox{\bf 5}}
\def\6{\mbox{\bf 6}}
\def\7{\mbox{\bf 7}}
\def\8{\mbox{\bf 8}}
\def\9{\mbox{\bf 9}}

\def\b{\mbox{\boldmath $b$}}
\def\cc{\mbox{\boldmath $c$}}

\def\e{\mbox{\boldmath $e$}}

\def\u{\mbox{\boldmath $u$}}
\def\v{\mbox{\boldmath $v$}}
\def\w{\mbox{\boldmath $w$}}
\def\x{\mbox{\boldmath $x$}}

\def\A{\mbox{\boldmath $A$}}
\def\B{\mbox{\boldmath $B$}}
\def\C{\mbox{\boldmath $C$}}

\def\F{\mbox{\boldmath $F$}}
\def\G{\mbox{\boldmath $G$}}
\def\H{\mbox{\boldmath $H$}}
\def\I{\mbox{\boldmath $I$}}

\def\O{\mbox{\boldmath $O$}}

\def\Q{\mbox{\boldmath $Q$}}

\def\W{\mbox{\boldmath $W$}}
\def\X{\mbox{\boldmath $X$}}
\def\Y{\mbox{\boldmath $Y$}}

\def\Inprod#1#2{\left\langle#1, \, #2\right\rangle}
\def\inprod#1#2{\langle#1, \, #2\rangle}

\def\Real{\mbox{$\mathbb{R}$}}
\def\coneK{\mbox{$\mathbb{K}$}}

\def\coneJ{\mbox{$\mathbb{J}$}}

\def\spaceV{\mbox{$\mathbb{V}$}}
\def\SymMat{\mbox{$\mathbb{S}$}}

\def\SymN{\mbox{$\mathbb{N}$}}

\def\DNN{\mbox{$\mathbb{D}$}}

\def\lb{\mbox{\rm lb}}
\def\lbb{\mbox{$\widetilde{\mbox{\rm lb}}$}}
\def\ub{\mbox{\rm ub}}

\begin{document}

\title{ \Large 
A Newton-bracketing method \\
for a simple conic optimization problem
}

\author{
\normalsize 
Sunyoung Kim\thanks{Department of Mathematics, Ewha W. University, 52 Ewhayeodae-gil, Sudaemoon-gu, Seoul 120-750, Korea 
			({\tt skim@ewha.ac.kr}). The research was supported
               by   NRF 2017-R1A2B2005119.}, \and \normalsize
Masakazu Kojima\thanks{Department of Industrial and Systems Engineering,
	Chuo University, Tokyo 192-0393, Japan ({\tt kojima@is.titech.ac.jp}).
	This research was supported by Grant-in-Aid for Scientific Research (A) 19H00808 %26242027
%	and the Japan Science and 
%	Technology Agency (JST), the Core Research of Evolutionary Science and 
%	Technology (CREST) research project.             
 	},
  \and \normalsize
Kim-Chuan Toh\thanks{Department of Mathematics, and Institute of Operations Research and Analytics, National University of Singapore,
10 Lower Kent Ridge Road, Singapore 119076
({\tt mattohkc@nus.edu.sg}). 
This research is supported in part by the Ministry of Education, Singapore, Academic Research Fund (Grant number: R-146-000-257-112).
} 
}

\date{\normalsize\today}

\maketitle 
\vspace*{-0.4cm}
\begin{abstract}
\noindent
%Long version:\\
For the Lagrangian-DNN relaxation of quadratic optimization problems (QOPs), we propose a Newton-bracketing method to improve the performance of the bisection-projection method implemented in BBCPOP [to appear in ACM Tran. Softw., 2019].  The relaxation problem is converted into the problem of finding the largest zero $y^*$ of a continuously differentiable (except at $y^*$) convex function $g : \Real \rightarrow \Real$ such that $g(y) = 0$ if $y \leq y^*$ and $g(y) > 0$ otherwise. In theory, the method generates lower and upper bounds of $y^*$ both converging to $y^*$. Their convergence is quadratic if the right derivative of $g$ at $y^*$ is positive.  Accurate computation of $g'(y)$ is necessary for the robustness of the method, but it is difficult to achieve in practice. As an alternative,  we present a secant-bracketing method. We demonstrate that the method improves the quality of the lower bounds obtained by BBCPOP and SDPNAL+ for binary QOP instances from BIQMAC. Moreover, new lower bounds for the unknown optimal values of large scale QAP instances from QAPLIB are reported. 
\end{abstract}

\noindent
{\bf Key words. } 
Nonconvex quadratic optimization problems,  conic relaxations, robust numerical algorithms,
Newton-bracketing method, secant-bracketing method for generating valid bounds.

\vspace{0.5cm}

\noindent
{\bf AMS Classification.} 
90C20,  	%Quadratic programming 
90C22,  	%Semidefinite programming
90C25, 	%Convex programming
%90C26.  	%Nonconvex programming, global optimization

%\input sect1.tex
\section{Introduction}

We consider a class of (nonconvex) quadratic optimization problems (QOPs) with linear and complementarity constraints in 
nonnegative continuous variables and binary variables. 
We are particularly interested in  efficient  and robust numerical algorithms for solving the class of QOPs.
This class includes many important combinatorial optimization 
problems such as binary integer QOPs, max-cut problems, maximum stable set problems, quadratic multi-knapsack problems 
and quadratic assignment problems. Solving a QOP in the class is known to be NP-hard
in general, and various (convex) conic optimization problem 
relaxations including popular semidefinite programming (SDP) relaxations 
have been studied by many researchers from both  theoretical and practical perspectives. 
%See \cite{ANJOS2012,FUJIE1997,GEOMANS95,VANDER96,WOLKOWICZ2000}. 
See \cite{ANJOS2012,GEOMANS95,SHOR1987,VANDER96,WOLKOWICZ2000}.%,SHOR1987}. 

One of notable developments in the theory of conic relaxations for QOPs is Burer's work in  \cite{BURER2009}, where
a completely positive programming (CPP) reformulation of a QOP 
in the class was proposed under an additional assumption on the feasible 
region of the QOP. This reformulation is very powerful in the sense that 
it attains the exact optimal value of the QOP. 
But unfortunately the CPP relaxation problem is numerically intractable.

Our focus in this paper is on tractable numerical methods. In particular, we are concerned with a Lagrangian doubly nonnegative (DNN) relaxation, which was proposed by 
Kim, Kojima and Toh  \cite{KIM2013} as a numerically tractable relaxation of the simplified Lagrangian-CPP reformulation \cite{ARIMA2012b} of a QOP in the class. 
The distinctive feature of their Lagrangian-DNN relaxation problem 
is the representation of the constraint set, which consists of 
a single linear equality constraint in a DNN variable matrix. 
They further reduced the optimal value of the dual 
of the relaxation problem to  the largest zero $y^*$ of an equation $g(y)=0$ 
in a single variable $y$. Here $g : \Real \rightarrow \Real$ 
satisfies the following properties: (a) $g$ is convex on $\Real$. 
(b) $g(y)=0$ if $y \leq y^*$ and $g(y) > 0$ otherwise. (c) 
Although $g$ is not explicitly represented, the function value 
$g(y)$ at each $y\in\Real$ can be computed through the metric projection onto the DNN cone. See Figure 1 (i). The properties (a) and (b) 
naturally lead to the bisection algorithm for approximately computing $y^*$. 
Step~0: Choose $\lb$ and $\ub$ such that $\lb < y^* < \ub$ and a sufficiently small positive 
$\epsilon > 0$. Step~1: Let $y = (\lb+\ub)/2$. 
Step~2: If $g(y) < \epsilon$ then $\lb = y$. Otherwise $\ub = y$. Go to Step 1. 
Based on this idea, they proposed the bisection-projection (BP) method 
for approximating the dual optimal value 
$y^*$ of the Lagrangian-DNN relaxation problem, and showed through 
numerical results that their method could efficiently compute  high quality lower bounds
for the optimal values of various QOP instances in the class.

\begin{figure}
\begin{center}
 \ifpdf
 \includegraphics[width=0.45\textwidth]{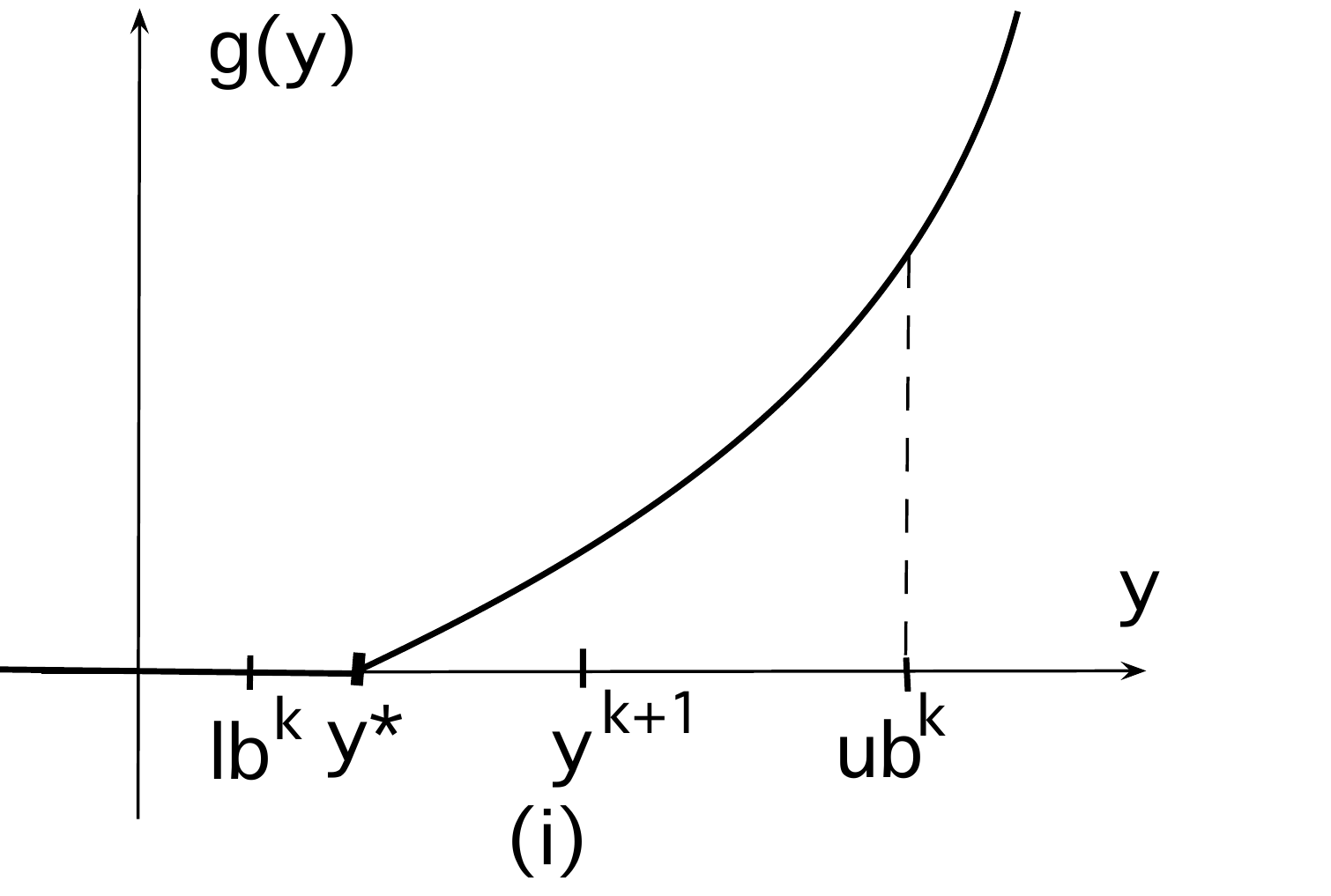}
 \else
 \includegraphics[width=0.45\textwidth]{bisection.eps}
 \fi
\ifpdf
\includegraphics[width=0.45\textwidth]{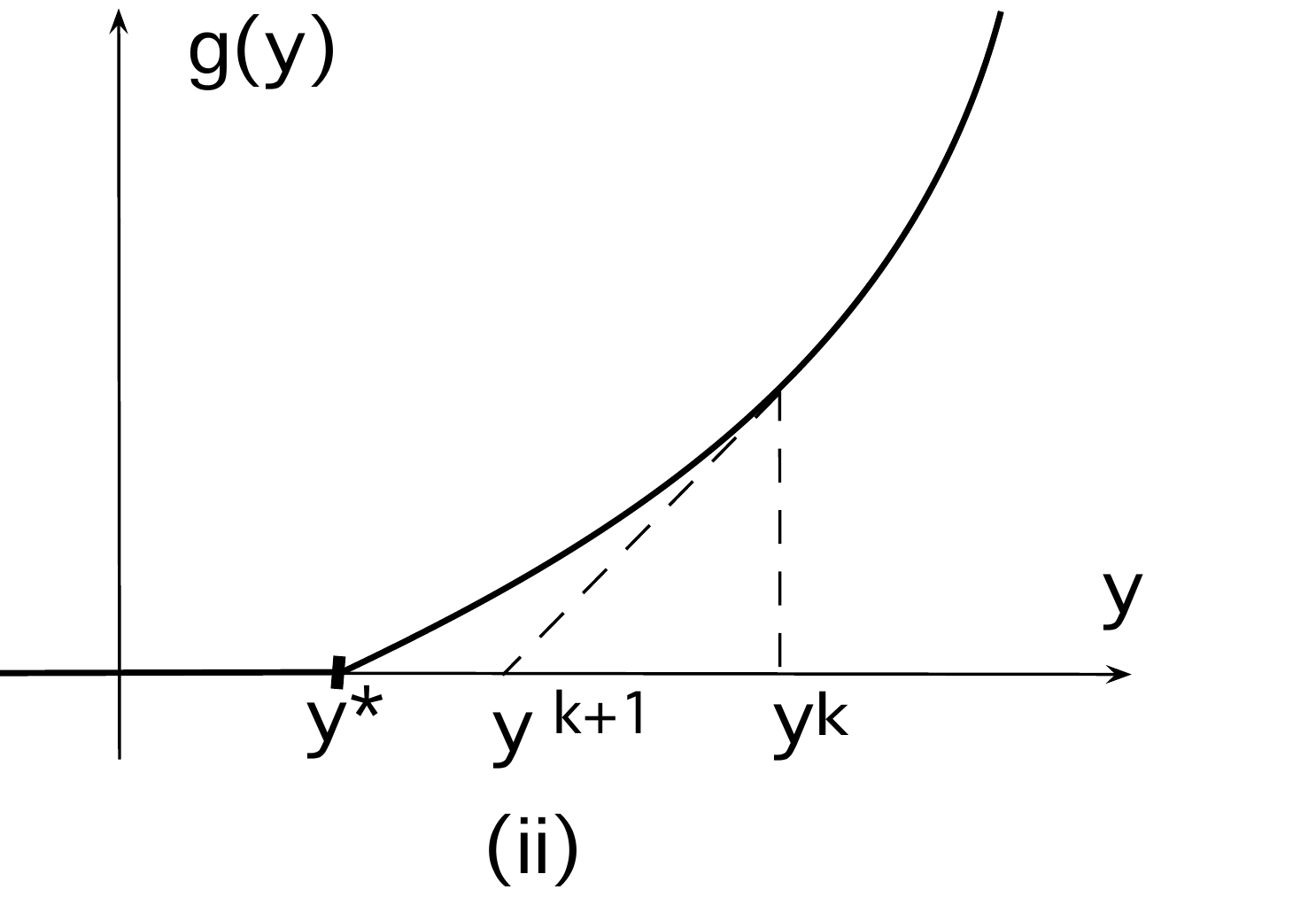}
\else
\includegraphics[width=0.45\textwidth]{Newton.eps}
\fi
\caption{
(i) The function $g : \Real \rightarrow \Real$ and the bisection algorithm.  
(ii) The Newton iteration.}
\label{figure:gy}
\end{center}
\end{figure}

Their simplified Lagrangian-DNN relaxation model and the BP method 
have been studied further to handle polynomial optimization problems (POPs) with binary, box and complementarity (BBC constraints), and extended to the Lagrangian conic 
optimization problem relaxation model in 
\cite{ARIMA2017,ARIMA2018,ARIMA2018b,KIM2016}, where the DNN cone is replaced 
with a more general convex cone. 
Recently, Ito, Kim, 
Kojima, Takeda and Toh \cite{ITO2018} released a software package 
BBCPOP for solving QOPs and POPs with BBC constraints. 
It was demonstrated in \cite{ITO2018} that BBCPOP can 
efficiently compute high quality lower bounds 
for optimal values of large scale QOPs and POPs with BBC constrains. 

The main purpose of this paper is to improve the performance 
of BBCPOP by incorporating the $1$-dimensional Newton method for approximating the 
largest zeros $y^*$ of the equation $g(y)=0$. See Figure 1 (ii). 
This idea was originally presented 
in \cite{ARIMA2018}, which showed that (d) $g : (y^*,\infty) \rightarrow \Real$ 
is continuously differentiable, in addition to (a), (b) and (c) mentioned above.  
However, neither the explicit algebraic representation 
of the derivative $g'(y)$ at $y \in \Real$  
nor computing its exact value is possible, although 
the derivative $g'(y)$ at $y \in \Real$ 
can be approximately computed as a by-product of the approximation of $g(y)$ through 
the metric projection onto the DNN cone. 
Due to this main difficulty, the $1$-dimensional Newton
was not incorporated in BBCPOP. 

Here we propose a Newton-bracketing method for the conic optimization problem (COP) 
converted from a given QOP
 using the technique proposed in  \cite{ARIMA2018,KIM2013}. More precisely,
we pay attention to the fact that the sequence $\{y^k\}$ generated by the 
Newton iteration $y^{k+1} = y^k - g(y^k)/g'(y^k)$ $(k=0,1,\ldots)$, starting from $y^0$ satisfying $g(y^0) > 0$, monotonically 
converges to $y^*$ from the above under the assumption that 
$g(y^k)$ and $g'(y^k)$ are exact. 
By itself, however, it does not generate any lower bound for $y^*$. By employing the technique 
proposed in \cite{ARIMA2017}, we generate a 
valid lower bound $\lbb$ at any iterate $y^k$ and the upper bound is updated by $y^{k+1}$.
Thus the Newton-bracketing method (Algorithm~\ref{algorithm:newton}) 
proposed in this paper generates both upper bound $y^k$ and lower bound $\lbb^k$ for $y^*$ at each iteration, 
and the upper and lower bounds converge to $y^*$ monotonically. Furthermore, we show that 
 their convergence is locally quadratic if the right derivative 
of $g$ at $y^*$ is positive (See Theorem~\ref{theorem:main}).

Despite the nice theoretical properties of 
the Newton-bracketing method  mentioned above,  its implementation 
 can be very expensive computationally in
practice since the accurate computation of the derivative $g'(y^k)$ is 
difficult. To avoid this difficulty, we present the damped 
secant-bracketing method (Algorithm~\ref{algorithm:secant}) by replacing the derivative by the secant 
$(g(y^k)-g(y^{k-1}))/(y^k-y{k-1})$. The proposed secant-bracketing method is implementable 
as the accurate computation of the function value 
is more manageable than the accurate evaluation of the derivative. 
We present some preliminary numerical results on the method 
applied to binary QOP instances from \cite{BIQMAC} and quadratic 
assignment problem (QAP) instances from \cite{QAPLIB}. We observe 
that the lower bounds obtained for the optimal values of the 
binary QOP instances by the new method are  tighter   
than those by BBCPOP and SDPNAL+. Also new  lower bounds obtained
by the method for large scale QAP instances (whose optimal values are not known)  are reported.

The paper is organized as follows. In Section 2.1, we introduce 
a primal-dual pair of conic optimization problems, 
COPs~\eqref{eq:COP1p} and~\eqref{eq:COP1d}, which serve as 
the aforementioned DNN relaxation of a QOP in the class and 
its dual. In Section 2.2, we present some fundamental properties (Lemma~\ref{Lemma4}) of  their Lagrangian relaxations described as COPs~\eqref{eq:COP1Lp} and~\eqref{eq:COP1Ld}. In Section 
2.3,  the BP method is described. We present 
the Newton-bracketing method in Section 3, the secant-bracketing 
method in Section 4, and  how a QOP in the class can be reduced to COPs~\eqref{eq:COP1p} and~\eqref{eq:COP1d}. 
Section 5. Some preliminary numerical results on the secant-bracketing method applied to binary QOP instances from \cite{BIQMAC} and large scale QAP instances from \cite{QAPLIB} are 
given in Section 6. We conclude in Section 7.

\section{Preliminaries}

\subsection{A primal-dual pair of simple conic optimization problems} 

Let $\Real$ denote the set of real numbers and $\Real_+$ the set of nonnegative real numbers. 
We use the following notation and symbols throughout the paper. 
\begin{eqnarray*}
\spaceV  & = & \mbox{a finite dimensional vector space endowed with an inner product } \nonumber \\ 
                 & \ & \mbox{$\inprod{\X}{\Y}$ and a norm $\|\X\| = \sqrt{\inprod{\X}{\X}}$ for every $\X, \ \Y \in \spaceV$},  \nonumber\\ [3pt]
\coneK_i    & = & \mbox{a nonempty closed convex cone in $\spaceV$} \ (i=1,2), \\ [3pt]
\coneK & = & \coneK_1 \cap \coneK_2. 
\end{eqnarray*}
For every nonempty closed convex cone $\coneJ$, 
its dual is denoted by $ \coneJ^*$, {\it i.e.}, 
$\coneJ^* = \left\{ \Y \in \spaceV : \right.$ 
$\left. \inprod{\X}{\Y} \geq 0 \right.$ $ \left. \mbox{for every } \X \in \coneJ \right\}$. 
In general, $\coneK^*$ coincides with the closure of $\coneK_1^* + \coneK_2^*$. Throughout the paper, we 
assume that $\coneK^* = \coneK_1^* + \coneK_2^*$. 

Let $\H^0, \H^1 \in \spaceV$. 
We consider the following pair of primal-dual COPs. 
\begin{eqnarray}
\eta^p & = & 
\inf \; \left\{ \inprod{\Q^0}{\X} : 
\begin{array}{l}
\X \in \coneK_1 \cap \coneK_2, \ \inprod{\H^0}{\X} = 1, \
\inprod{\H^1}{\X} = 0 
\end{array}
\right\}, \label{eq:COP1p} \\ 
\eta^d & = & 
\sup \; \left\{ y_0 : 
\begin{array}{l}
\displaystyle 
\Q^0 - \H^0 y_0  + 
\H^1y_1 = \Y_1 + \Y_2, \ \Y_1 \in \coneK_1^*, \ \Y_2 \in \coneK_2^*
\end{array}
\right\}.   \label{eq:COP1d}
\end{eqnarray}
We assume the following condition.
\begin{description}
\item[Condition (I) ] 
COP~\eqref{eq:COP1p} is feasible, 
and $\H^0, \H^1 \in \coneK_1^* + \coneK_2^*$. 
\end{description}

\subsection{A Lagrangian relaxation of COP~\eqref{eq:COP1p} and its dual}

If the Lagrangian relaxation is applied to COP~\eqref{eq:COP1p}, 
the following COP and its dual are obtained.
\begin{eqnarray}
\eta^p_{\lambda} & = & 
\inf \; \left\{ \inprod{\Q^0 + \lambda \H^1}{\X} :  
\begin{array}{l}
\X \in \coneK_1 \cap \coneK_2, \ \inprod{\H^0}{\X} = 1
\end{array}
\right\}, \label{eq:COP1Lp} 
\\[3pt] 
\eta^d_{\lambda} & = & 
\sup \; \left\{ y_0 : 
\begin{array}{l}
\displaystyle 
\Q^0 + \lambda \H^1 - \H^0 y_0 = \Y_1 + \Y_2, \ \Y_1 \in \coneK_1^*, \ \Y_2 \in \coneK_2^*
\end{array}
\right\},  \label{eq:COP1Ld}
\end{eqnarray}
where $\lambda \in \Real$ denotes the Lagrangian multiplier for the homogeneous equality $\inprod{\H^1}{\X} = 0$ in \eqref{eq:COP1p}. Since $\H^1 \in \coneK_1^*+\coneK_2^*$ by
Condition (I), the term $\lambda\inprod{\H^1}{\X}$ added to the 
objective function is nonnegative for every $\X \in \coneK_1\cap\coneK_2$ 
and $\lambda \geq 0$. As a result, the term serves as a penalty term for 
the violation of $\inprod{\H^1}{\X} = 0$ such that if 
$\X \in \coneK_1\cap\coneK_2$ and $\inprod{\H^1}{\X} \not= 0$ then 
$\lambda\inprod{\H^1}{\X} \rightarrow \infty$ as 
$\lambda \rightarrow \infty$. Furthermore, the following 
lemma holds. 

\lemm \cite[Lemmas 2.3 and 2.5]{ARIMA2018}
\label{Lemma4} Suppose that Conditions (I) is satisfied. 
 Then, the following assertions hold. 
\begin{description}
\item{(i) } $\eta^d_{\lambda} = \eta^p_{\lambda}$ for every $\lambda \in \Real$. 
Moreover, if 
 $\eta^p_{\lambda}$ is finite, then \eqref{eq:COP1Ld} has an optimal solution with the 
objective value $\eta^d_{\lambda} = \eta^p_{\lambda}$. 
\item{(ii) } $\displaystyle
 \left( \eta^d_{\lambda} = \eta^p_{\lambda} \right) \hspace{-1.5mm}\uparrow  
= \eta^d$. Here $\uparrow  = \eta^d$ means ``increases monotonically and converges to $\eta^d$ as $\lambda \rightarrow \infty$''. 
\item{(iii) } Assume in addition that $\{ \X \in F : \inprod{\Q^0}{\X} \leq \bar{\eta} \}$ is bounded for some 
$\bar{\eta}$. Then $\eta^d = \eta^p$. 
\end{description}
\elemm

In the remaining of the paper, we impose the following condition 
on the space $\spaceV$ and the cone $\coneK=\coneK_1\cap\coneK_2$.
\vspace{-2mm}
\begin{description}
\item[Condition (II) ]
$\spaceV$ is the Cartesian product of symmetric matrix spaces, 
 % $\coneK =\coneK_1 \cap \coneK_2 \subset \spaceV$ with 
 $\coneK_1$ the cone 
of consisting of the positive semidefinite matrices in $\spaceV$,  and $\coneK_2$ a 
closed convex cone in $\spaceV$ such that $(\coneK_1 \cap \coneK_2)^* = \coneK_1^* + \coneK_2^*$.\vspace{-2mm}
\end{description}
In this case, we know that $\coneK_1^* = \coneK_1$. 
By Lemma \ref{Lemma4}, we can obtain an accurate lower bound 
$\eta^d_\lambda$ for the optimal value $\eta^d$ of COP~\eqref{eq:COP1d} 
by solving COP~\eqref{eq:COP1Ld} for 
a sufficiently large positive $\lambda$. 
For simplicity of notation, we fix $\lambda$ to be a sufficiently 
large positive number, and we rewrite the primal-dual pair of 
COPs~\eqref{eq:COP1Lp} and~\eqref{eq:COP1Ld} as 
\begin{eqnarray}
\varphi^* & = & 
\inf \; \left\{ \inprod{\Q}{\X} :  
\begin{array}{l}
\X \in \coneK_1\cap\coneK_2, \ \inprod{\H}{\X} = 1
\end{array}
\right\}, \label{eq:COP2Lp} 
\\[3pt] 
y^* & = & 
\sup \; \left\{ y : 
\begin{array}{l}
\displaystyle 
\Q - \H y  = \Y_1 + \Y_2, \ \Y_1 \in \coneK_1^*, \ \Y_2 \in \coneK_2^*
\end{array}
\right\}.  \label{eq:COP2Ld}
\end{eqnarray}

To solve COP~\eqref{eq:COP2Ld}, we  present the bisection-projection (BP) method (Algorithm~\ref{algorithm:bisection}) 
 in Section 2.3, 
the Newton-bracketing method (Algorithm~\ref{algorithm:newton}) in Section 3 
and the secant-bracketing method (Algorithm~\ref{algorithm:secant}) in Section 4. 

%-----------------------------------------------------------------------------------------------------------------
\subsection{The bisection-projection method for COP~\eqref{eq:COP2Ld}}

We describe the bisection-projection (BP) method 
\cite{ARIMA2017,ARIMA2018,ITO2018,KIM2013,KIM2016} 
for COP \eqref{eq:COP2Ld}. 
For every $y \in \Real$, define $\G(y) = \Q - \H y$ and
\begin{eqnarray}
% \G(y) & = & \Q - \H y, \label{eq:defG} \nonumber \\[3pt] 
g(y) & = & \min \; \left\{ \left\| \G(y) - (\Y_1+\Y_2) \right\| :  \Y_1\in \coneK_1^*, \ \Y_2 \in \coneK_2^* \right\} 
\label{eq:defgy0} \\
& = &  \| \G(y) - (\widehat{\Y}_1(y)+\widehat{\Y}_2(y)) \|. \nonumber
\end{eqnarray}
Here $(\widehat{\Y}_1(y),\widehat{\Y}_2(y)) \in \coneK_1^* \times \coneK_2^*$ denotes an optimal solution. 
For every $y \in \Real$, we have that $g(y) \geq 0$, and 
\begin{eqnarray*}
g(y)=0 & \Leftrightarrow & \G(y) = \Q -\H y \in \coneK_1^* + \coneK_2^* \\ 
& \Leftrightarrow & \mbox{$y$ is a feasible solution of COP~\eqref{eq:COP2Ld}}.
\end{eqnarray*}
Thus,  COP~\eqref{eq:COP2Ld} can be rewritten as 
\begin{eqnarray}
y^* & = & \sup \; \left\{ y \in \Real : g(y) = 0  \right\}. \label{eq:equalityg}
\end{eqnarray}
From $\H \in \coneK_1^*+\coneK_2^*$ by Condition (I), we see that $g(y)=0$ ({\it i.e.}, $\Q -\H y \in \coneK_1^*+\coneK_2^*$)
if $y \leq y^*$ and $g(y)>0$ ({\it i.e.},  $\Q -\H y \not\in \coneK_1^*+\coneK_2^*$) otherwise. 

In addition to Conditions (I) and (II), we assume the following condition in the subsequent discussion.\vspace{-3mm}
\begin{description}
\item[Condition (III) ] (Condition (D) of \cite{ARIMA2017}) A positive number $\rho$ such that $ \inprod{\I}{\X}  \leq \rho$ for all feasible solutions $\X$ of \eqref{eq:COP2Ld} is 
known, where $\I$ denotes the identity matrix in $\spaceV$.\vspace{-3mm}
\end{description}
Then, the primal-dual pair of COPs \eqref{eq:COP2Lp} and \eqref{eq:COP2Ld} 
are equivalent to 
\begin{eqnarray}
& & 
\begin{array}{llll}
\varphi^* = \min \left\{ \inprod{\Q}{\X} : \inprod{\H}{\X} = 1, \ \Inprod{\I}{\X} \leq \rho, \  \X \in \coneK_1 \cap \coneK_2 \right\}, 
\end{array}
\label{eq:COP2}
\end{eqnarray} 
and its dual 
\begin{eqnarray}
& & \begin{array}{llll}
y^* = \max \left\{ y +\rho t : \G(y) -\I t -\Y_2 = \Y_1 \in \coneK_1^*, \ \Y_2 \in \coneK_2^*, \ t \leq 0 
\right\}, 
\end{array}
\label{eq:COPD2}
\end{eqnarray}
respectively. For every $y \in \Real$ and $\Y_2 \in \coneK_2^*$, let 
\begin{eqnarray}
\left\{
\begin{array}{lcl}
\tilde{t}(y,\Y_2) =  \min\{0, \ \lambda_{\min}(\G(y) - \Y_2)\}, \\
\widetilde{\Y}_1 (y,\Y_2) = \G(y) -\I \tilde{t}(y,\Y_2) -\Y_2, 
\end{array}
\right.
\label{eq:tildet}
\end{eqnarray}
where $\lambda_{\min}(\A)$ denotes the minimum eigenvalue of $\A \in \spaceV$. 

\lemm \label{lemma:LB1} (Lemma 3.1 of \cite{ARIMA2017})
Let $y \in \Real$ and $\Y_2 \in \coneK_2^*$. Then, $y + \rho\tilde{t}(y,\Y_2) \leq y^*$. 
\elemm
\proof{
By the construction \eqref{eq:tildet} of $\tilde{t}(y,\Y_2)$ and $\widetilde{\Y}_1 (y,\Y_2)$, we see that 
\begin{eqnarray*}
\G(y) -\I \tilde{t}(y,\Y_2) -\Y_2 = 
\widetilde{\Y}_1 (y,\Y_2) \in \coneK_1^*, 
\end{eqnarray*}
which, together with $\Y_2 \in \coneK_2^*$, implies that 
$(t,y,\Y_1,\Y_2)$ with $ t = \tilde{t}(y,\Y_2)$ and $\Y_1 = \widetilde{\Y}_1 (y,\Y_2)$ 
is a feasible solution of \eqref{eq:COPD2} with the objective value $y+\rho \tilde{t}(y,\Y_2)$. 
Hence, the inequality $y+\rho \tilde{t}(y,\Y_2) \leq y^*$ follows.
\qed
}

In Algorithm~\ref{algorithm:bisection} described below, we utilize Lemma \ref{lemma:LB1} to generate 
a valid lower bound $\lbb$ for $y^*$ under the assumption that the computation of the minimal eigenvalue of 
$\lambda_{\min}(\A)$ is accurate for every $\A \in \spaceV$. 

Let $y \in \Real$. Then, every optimal solution 
$(\Y_1,\Y_2) = (\widehat{\Y}_1(y),\widehat{\Y}_2(y))$ 
of the minimization problem~\eqref{eq:defgy0} is characterized by the KKT 
condition: there exists an $\X = \widehat{\X}(y)$ such that  
\begin{eqnarray}
& & \G(y) = \Y_1 + \Y_2- \X,  
\ \Y_1 \in \coneK_1^*, \ \Y_2 \in \coneK_2^*, \label{eq:KKT0}  \\ 
& & \X \in \coneK_1 \cap \coneK_2, \ \inprod{\X}{\Y_1} = 0  \   \mbox{and }
\inprod{\X}{\Y_2} = 0. 
\label{eq:KKT1} 
\end{eqnarray}
Therefore, we obtain 
\begin{eqnarray}
g(y) & = & \| \G(y) - (\widehat{\Y}_1(y) + \widehat{\Y}_2(y))\| = 
\| \widehat{\X}(y)  \|. \label{eq:gy1}
\end{eqnarray}

\algo \label{algorithm:bisection} 
(the BP method \cite[Algorithm 3.2]{ARIMA2017})\vspace{-2mm}
\begin{description}
\item{Step 0: } Let $\delta$ and $\epsilon$ be sufficiently small positive numbers. 
Here $\delta > 0$ 
determines the length of the target interval  $[\lb,\ub]$ 
and $\epsilon > 0$  
is used to determines whether $g(y)$ attains $0$ numerically 
with $\|\X\| < \epsilon$ in Step 3. See\eqref{eq:gy1}. For instance,
we can use $\delta = 1.0$e-4 and $\epsilon = 1.0$e-12 in the double precision arithmetic. 
Choose $\lbb$, $\lb$ and  $\ub$ such that 
$\lbb = \lb < y^* < \ub$, where $\lbb =\lb$ can be $-\infty$. \vspace{-2mm}
\item{Step 1:  } If $\ub - \lb < \delta$, then output $\lbb$ as a valid lower bound for 
$y^*$ and $[\lb,\ub]$ as an interval which expects to contain $y^*$. \vspace{-2mm}
\item{Step 2:  } Let $y = (\lb+\ub)/2$.\vspace{-2mm}
\item{Step 3:  } 
Compute $(\X,\Y_1,\Y_2)$ which satisfies accurately~\eqref{eq:KKT0} 
 and approximately ~\eqref{eq:KKT1} by applying the accelerated proximal gradient (APG) method \cite{BECK2009}  to the problem~\eqref{eq:defgy0}. 
Let % $g(y) = \|\X\|$ and 
$\lbb = \max\{\lbb , $ $y + \rho\tilde{t}(y,\Y_2)\}$.  
Let $\lb = y$ if $\|\X\| < \epsilon$, and $\ub = y$ otherwise ({\it i.e.}, $\|\X\| \geq \epsilon$). Go to Step 1.
\end{description}
\ealgo

At Step 3,  $\Y_2 \in \coneK_2^*$ is obtained. As a result, 
Lemma~\ref{lemma:LB1} guarantees that $\lbb$  is 
a valid lower bound for $y^*$ 
regardless of the choice of a positive $\epsilon$. On the other hand, 
$y^* \in [\lb,\ub]$ is not guaranteed because of the numerical error.

%\red{
%We can compute 
%$(\X,\Y_1,\Y_2)$ which satisfies~\eqref{eq:KKT0} 
%accurately and~\eqref{eq:KKT1} approximately} by applying 
The APG method \cite{BECK2009} used in Step 3 is a first order iteration method  
that employs the metric projection $\Pi_i^*$ onto $\coneK_i^*$ $(i=1,2)$.
% , to the problem~\eqref{eq:defgy0}. 
The BP method combined with 
the APG method has been  implemented recently as a software package 
BBCPOP \cite{ITO2017}. 
It employs the following error criterion % based on the KKT condition ~\eqref{eq:KKT0}. 
\begin{eqnarray}
\kappa(\X,\Y_1,\Y_2) &=& \max \left\{
\displaystyle 
\begin{array}{l}
\displaystyle
\frac{\inprod{\X}{\Y_1}}{1+\|\X\|+\|\Y_1\|}, \
\frac{\inprod{\X}{\Y_2}}{1+\|\X\|+\|\Y_2\|}, \vspace{2mm}\\ 
\displaystyle
\frac{\|\Pi_1^*(-\X)\|}{1+\|\X\|}, \
\frac{\|\Pi_2^*(-\X)\|}{1+\|\X\|}
\end{array}
\right\}
\end{eqnarray}
to decide whether the iterate 
$(\X,\Y_1,\Y_2)$ of the APG method  approximately satisfies% ~\eqref{eq:KKT0} accurstely and
~\eqref{eq:KKT1}. 
We note that~\eqref{eq:KKT0} is 
maintained 
throughout the iterations, and that $\Pi_i^*(-\X) = \X - \Pi_i(\X)$ $(i=1,2)$ holds for every $\X \in \spaceV$ by Moreau's decomposition theorem \cite{MOREAU62}. 
Here $\Pi_i(\X)$ ($\Pi_i^*(\X)$, respectively) denotes the metric projection of $\X\in\spaceV$ 
onto $\coneK_i$ ($\coneK_i^*$, respectively) $(i=1,2)$. If $\|\X\| < \epsilon$ or $\|\X\| \geq \epsilon$ 
and $\kappa(\X,\Y_1,\Y_2) < tol$, then the iteration stops, where 
$\epsilon$ and 
$tol$ are sufficiently small positive numbers, such as 
$\epsilon = tol = 1.0$e-12 in the double precision arithmetic. 
See Sections of~2.2 and~4.2 of \cite{ITO2018} for more details.

In the APG method, the metric projections from $\spaceV$ onto 
the cones $\coneK_i$ $(i=1,2)$  and their duals play an essential role. 
We have been dealing with the case where 
$\coneK_1 = \coneK_1^*$ is the cone of symmetric matrices in 
$\spaceV$, thus  the metric projection of $\A \in \spaceV$ onto  $\coneK_1 = \coneK_1^*$ 
is carried out efficiently and accurately 
via the eigenvalue decomposition of $\A$. For the efficient 
and accurate metric projection onto the cone $\coneK_2$, 
we need to restrict $\coneK_2$ to a class of polyhedral cones 
onto which the accurate metric projection can be efficiently computed. Such a class of polyhedral cones  were studied in \cite{KIM2016}. The class includes the cone of nonnegative matrices 
in $\spaceV$, polyhedral cones induced from SDP relaxations of binary and 
complementarity constraints, and their intersections. Those cones are described in 
Section~5 where applications of COP~\eqref{eq:COP1p} to doubly nonnegative 
relaxations of combinatorial QOPs are discussed.

\section{The Newton-bracketing method}

We use the same notation and symbols as in Section 2, specifically 
$(\widehat{\Y}_1(y),\widehat{\Y}_2(y))$ denotes an optimal solution of the minimization 
problem~\eqref{eq:defgy0} and  
$(\widehat{\X}(y),\widehat{\Y}_1(y),\widehat{\Y}_2(y)) \in (\coneK_1\cap\coneK_2) \times \coneK_1^* \times \coneK_2^*$ satisfies the KKT condition~\eqref{eq:KKT0} and~\eqref{eq:KKT1}; 
See also~\eqref{eq:gy1}. 
To describe the $1$-dimensional Newton method for computing 
$y^*$, 
we need the following lemma, 
which shows some fundamental properties of the function $g : \Real \rightarrow \Real_+$. 
\lemm ( \cite[Lemma 4.1]{ARIMA2018}) \label{lemma:y0} 
Assume that $y^*$ is finite. \vspace{-2.0mm}
\begin{description} 
\item{(i) } $g : \Real \rightarrow \Real_+$ is continuous and convex. \vspace{-2.0mm} 
\item(ii)
If $ y > y^*$, then $\inprod{\H}{\widehat{\X}(y)} > 0$ and $\widehat{\X}(y) / \inprod{\H}{\widehat{\X}(y)} $
 is a feasible solution of the primal COP~\eqref{eq:COP2Lp}. \vspace{-2.0mm}
\item{(iii) } 
If $y > y^*$, then 
$g'(y) = dg(y) / dy = \inprod{\H}{\widehat{\X}(y)}/g(y) > 0$; 
hence $g : (y^*,\infty) \rightarrow \Real$ is continuously differentiable and strictly increasing. \vspace{-2.0mm}
\item{(iv) } 
Assume that $\G(\bar{z})$ lies in the interior of $\coneK_1^* \times \coneK_2^*$ for some $\bar{z}$. Then 
$\displaystyle \frac{g(y) - g(y^*)}
{y - y^*}$ converges to a positive value as 
$y \downarrow y^*$; hence the right 
derivative $g'_+(y^*)$ of $g(y)$ at $y = y^*$ is positive. \vspace{-2.0mm}
\end{description}
\elemm

Suppose that $g(\bar{y})>0$, {\it i.e.}, $\bar{y} > y^*$ for 
some $\bar{y} \in \Real $. Then the Newton iteration for computing 
$y^*$ is given by 
\begin{eqnarray*}
\bar{y}^+ & = & \bar{y} - \frac{g(\bar{y})}{g'(\bar{y})} 
\;=\; \bar{y} - \frac{\inprod{\widehat{\X}(\bar{y})}
{\widehat{\X}(\bar{y})}}{\inprod{\H}{\widehat{\X}(\bar{y})} } 
\\ 
& = & \bar{y} 
- \frac{\inprod{\widehat{\Y}_1(\bar{y}+\widehat{\Y}_2(y)) - \G(\bar{y})}
{\widehat{\X}(\bar{y})}}
{\inprod{\H}{\widehat{\X}(\bar{y})} } \quad \mbox{(since $\widehat{\X}(\bar{y}) 
= \widehat{\Y}_1(\bar{y})+ \widehat{\Y}_2(\bar{y})- \G(\bar{y})$ by \eqref{eq:KKT0})} \\ 
& = & \bar{y} + \frac{\inprod{\G(\bar{y})}
{\widehat{\X}(\bar{y})}}{\inprod{\H}{\widehat{\X}(\bar{y})} } \quad 
\mbox{(since $\inprod{\widehat{\X}(\bar{y})}{\widehat{\Y}_1(\bar{y})} = \inprod{\widehat{\X}(\bar{y})}{\widehat{\Y}_2(\bar{y})} = 0 $ by \eqref{eq:KKT1})} \\ 
& = &  \bar{y} + \frac{\inprod{\Q - \H\bar{y}}
{\widehat{\X}(\bar{y})}}{\inprod{\H}{\widehat{\X}(\bar{y})} } \quad 
\mbox{(by the definition %\eqref{eq:defG} 
                         of $\G(\bar{y}$))} \\ 
& = & \Inprod{\Q}{ \widetilde{\X}(\bar{y})} \;\geq\; \varphi^*,
\end{eqnarray*}
where $\widetilde{\X}(\bar{y}) = {\widehat{\X}(\bar{y})}/{\inprod{\H}{\widehat{\X}(\bar{y})}}$ 
denotes a feasible solution of the primal COP \eqref{eq:COP2Lp}. 
It should be noted that the sequence $\{ y^k \}$ generated by the Newton 
iteration  from any initial iterate $y^0 > y^*$ monotonically decreases and converges to $y^*$ 
by (i), (ii) and (iii) of Lemma~\ref{lemma:y0}. 

Now we combine the Newton iteration with Lemma~\ref{lemma:LB1} 
to design %for the 
Newton-bracketing method for solving COP~\eqref{eq:COP2Ld}.
\algo \label{algorithm:newton} \mbox{\ }\vspace{-2mm}
\begin{description}
\item{Step 0: }  Let $\lbb^0 = -\infty$ and $k = 0$. 
Choose a $y^0 > y^*$.\vspace{-2mm}
\item{Step 1: } Compute 
$(\X,\Y_1,\Y_2) % = (\widehat{\X}(y^k),\widehat{\Y}_1(y^k),\widehat{\Y}_2(y^k)
\in (\coneK_1\cap\coneK_2) \times \coneK_1^* \times \coneK_2^*$ which satisfies 
the KKT condition~\eqref{eq:KKT0} and~\eqref{eq:KKT1} with $y=y^k$; 
hence~\eqref{eq:gy1} holds with $(\widehat{\X}(y),\widehat{\Y}_1(y),\widehat{\Y}_2(y)) = (\X,\Y_1,\Y_2)$ and 
$y=y^k$. 
Let $g(y^k)=\|\X\|$.\vspace{-2mm}
\item{Step 2: } (Application of Lemma~\ref{lemma:LB1} to update $\lbb^k$) 
Let $\lbb^{k+1} = \max\{\lbb^{k}, \  y^k +\rho \tilde{t}(y^k,\Y_2 \}$.\vspace{-2mm}
\item{Step 3: } 
If $g(y^k) = 0$ then $y^k = y^*$ and stop the iteration.
\vspace{-2mm}
\item{Step 4: } (Newton step to update $y^k$) Let $g'(y^k) = \inprod{\H}{\X}/g(y^k)$, and 
$y^{k+1} = y^k - g(y^k)/g'(y^k)$.\vspace{-2mm}
\item{Step 5: } Replace $k$ by $k+1$, and go to Step 1.% \vspace{-2mm}
\end{description}
\ealgo 

\theo \label{theorem:main} \mbox{\ }
\vspace{-2mm}
\begin{description}
\item{(i) } $y^* \in [\lbb^k,y^k]$ for every $k=1,2,\ldots$, $y^k\downarrow = y^*$ 
and $\lbb^k\uparrow=y^*$. Here $\uparrow  = y^*$ (or $\downarrow  = y^*$) means ``increases (or decreases) monotonically and converges to $y^*$ as $k \rightarrow \infty$''
%and both $\lbb^k$ and 
%$y^k$ converge to $y^*$ monotonically as $k \rightarrow \infty$;
%$y^k > y^{k+1} \rightarrow y^*$ and  
%$\lbb^k \leq \lbb^{k+1} \rightarrow y^*$ as $k \rightarrow \infty$. 
\vspace{-2mm}
\item{(ii) }  Assume that $\G(\bar{z}_0)$ lies in the interior of $\coneK_1^*+\coneK_2^*)$ for some 
$\bar{z}_0$. Then the convergence of $y^k$ to $y^*$
as well as the convergence of $\lbb^k$ to $y^*$ are quadratic.
% \vspace{-2mm}
\end{description} 
\etheo

To prove the theorem, we need the following lemma.  
\lemm \label{lemma:LB2} \mbox{ \ } 
Assume that $g(y) > 0$ for some $y \in \Real$. Then 
$y - \rho g(y) \leq y + \rho\tilde{t}(y,\widehat{\Y}_2(y)) \leq y^*$. 
\elemm
\proof{
By Lemma~\ref{lemma:LB1}, it suffices to show that $-g(y) \leq \tilde{t}(y,\widehat{\Y}_2(y))$. 
It follows from \eqref{eq:KKT0} that 
\begin{eqnarray*}
\lambda_{\min}(\G(y) - \widehat{\Y}_2(y)) & = & \lambda_{\min}(\widehat{\Y}_1(y) - \widehat{\X}(y)) \\ 
& = & \u^T (\widehat{\Y}_1(y) - \widehat{\X}(y)) \u \ \mbox{for some $\u$ with $\left\|\u\right\| = 1$} \\
& \geq & -\u^T \widehat{\X}(y)) \u \ \mbox{(since $\widehat{\Y}_1(y) \in \coneK_1^*$)} \\
& \geq & - \lambda_{\max}(\widehat{\X}(y)) \\ 
& \geq & -\left\| \widehat{\X}(y) \right\| = - g(y). 
\end{eqnarray*}
Since $- g(y) < 0 $, 
$- g(y)\leq \min\{0, \ \lambda_{\min}(\G(y)-\widehat{\Y}_2(y))\} =
\tilde{t}(y,\widehat{\Y_2}(y))$ 
follows. 
\qed
}

\medskip 

\noindent
%\proofof{Theorem~\ref{theorem:main}}{
{\em Proof of Theorem~\ref{theorem:main}}. 
The assertions of (i) follows from Lemma~\ref{lemma:LB1}, Lemma~\ref{lemma:y0}, and the construction of the sequences $\{\lbb^k\}$ and 
$\{y^k\}$ by Algorithm~\ref{algorithm:newton}. To prove (ii), assume that 
$\G(\bar{z}_0)$ lies in the interior of $\coneK_1^*+\coneK_2^*$ for some 
$\bar{z}_0$. By (iv) of Lemma~\ref{lemma:y0}, the right 
derivative $g'_+(y^*)$ of $g(y)$ at $y = y^*$ is positive. 
Define 
\begin{eqnarray*}
h(y) & = & \left\{\begin{array}{lll}
g(y) & \mbox{if } y \geq y^*, \\
g'_+(y^*)(y-y^*) & \mbox{otherwise}.
\end{array}\right.
\end{eqnarray*}
Then $h : \Real \rightarrow \Real$ is continuously differentiable and convex, 
and $h'(y^*) = g'_+(y^*) > 0$. Since the application of Newton method 
to $h(y) = 0$ with 
the initial point $y^0$ satisfying $h(y^0) = g(y^0) > 0$ 
yields the same sequence as the Newton method applied to solve $g(y) = 0$, $y^k$ converges  quadratically  to 
$y^*$  (hence $g(y^k)$ converges to $0$)  as $k \rightarrow \infty$ by 
\cite[Monotone Newton Theorem 13.3.4]{ORTEGA1970}. 
This implies the quadratic convergence of $\{g(y^k)\}$ to $0$. By 
%Thus the desired result follows from 
Lemma~\ref{lemma:LB2}, $\lbb^k$ also converges quadratically to $y^*$. % quadratically too. 
\qed
%}

\section{The damped secant-bracketing method}

Although Algorithm~\ref{algorithm:newton} (the Netwton-bracketing method) 
has nice theoretical properties %in theory 
as shown in Theorem~\ref{theorem:main}, 
its % \red{stable and efficient} 
accurate implementation is difficult. %  \red{in practice}. 
More precisely, 
 the exact function value 
$g(y^k)$ and its derivative $g'(y^k)$ have been assumed for each iterate $y^k$ in Algorithm~\ref{algorithm:newton}. 
However, when Algorithm~\ref{algorithm:newton} is implemented in practice, only approximate value $(\X,\Y_1,\Y_2)$ 
of $(\widehat{\X}(y^k),\widehat{\Y}_1(y^k),\widehat{\Y}_2(y^k))$ can be computed by the APG method. 
From those values,  $g(y^k)$ and $g'(y^k)$ are approximated by 
$\|\X\|$ and $\inprod{\H}{\X}/\|\X\|$, respectively.
In particular, the accurate computation of $g'(y^k) = \inprod{\H}{\widehat{\X}(y^k)}/\|\widehat{\X}(y^k)\|$ is 
difficult as the denominator $\|\widehat{\X}(y^k)\|$ as well as the numerator $\inprod{\H}{\widehat{\X}(y^k)}$ both converge to zero 
as $y^k$ approaches to $y^*$. This computational issue should be 
carefully dealt with for the robustness
of the Newton-bracketing method in practice.

To avoid such a difficulty and for the purpose of  developing a practically robust method, 
we replace the derivative $g'(y^k)$ at the $k$th iterate $y^k$ 
by the secant $(g(y^k)-g(y^{k-1}))/(y^k-y^{k-1})$ $(k=1,2,\ldots)$. Then the iterative formula turns out to be 
\begin{eqnarray}
y^{k+1} & = & y^k -\alpha^k \frac{g(y^k)(y^k-y^{k-1})}{g(y^k)-g(y^{k-1})} 
\nonumber \\ 
& = &  y^k -\alpha^k \frac{\|\widehat{\X}(y^k)\|(y^k-y^{k-1})}
{\|\widehat{\X}(y^k)\|-\|\widehat{\X}(y^{k-1})\|}  \ (k=1,2,\ldots). \label{eq:secant}
\end{eqnarray}
Here $\alpha^k \in (0,1]$ is a damping factor, which is multiplied to the secant 
to avoid the occurrence of $y^{k+1} < y^*$  in any case, 
which might happen due to inaccurate computations of $g(y^{k-1})$ and $g(y^k)$. 
To initiate the first secant, we need to prepare two points $y^0 > y^1 > y^*$. 
% \green{Let $\X^k = \X(y^k)$}.

\algo \label{algorithm:secant} ({\bf The damped secant-bracketing method}) \vspace{-2mm}
\begin{description}
\item{Step 0: } Let $\delta$ and 
$\epsilon$ be sufficiently small positive numbers (for example, $\delta = 1.0$e-4 and 
$\epsilon = 1.0$e-12 in the double precision arithmetic). 
Choose $y^0$ and $y^1$  such that $y^* < y^1 < y^0$.  Let $\lbb^0 = -\infty$ 
and $k=0$. 
\vspace{-2mm}
\item{Step 1: } Compute 
$(\X,\Y_1,\Y_2) = (\X^k,\Y^k_1,\Y^k_2) $ 
which satisfies accurately ~\eqref{eq:KKT0}  and approximately~\eqref{eq:KKT1} 
 with $y=y^k$ by the APG method. 
% Let $g(y^k)=\|\widehat{\X}(y^k)\|$.
\vspace{-2mm}
\item{Step 2: } (Application of Lemma~\ref{lemma:LB1} to update $\lbb^k$) 
Let $\lbb^{k+1} = \max\{\lbb^{k}, \  y^k +\rho \tilde{t}(y^k,\Y^k_2\}$. 
If k=0, then let k=1 and go to Step1.\vspace{-2mm}
\item{Step 3: } If $\|\X^k\| < \epsilon$ or $y^k - \lbb^{k+1}  < \delta$, then 
stop the iteration, and output $\lbb^{k+1}$ as a valid lower bound for $y^*$ and $y^k$ as an 
approximation of $y^*$.\vspace{-2mm}
\item{Step 4: } (Damped secant iteration) 
Apply 
\begin{eqnarray*}
y^{k+1} & = & y^k -\alpha^k \frac{\|\X^k\|(y^k-y^{k-1})}{\|\X^k\|-\|\X^{k-1}\|} 
\end{eqnarray*}
with some $\alpha^k \in (0,1]$ to update $y^k$. \vspace{-2mm}
\item{Step 5: } Replace $k$ by $k+1$, and go to Step 1.  % \vspace{-2mm}
\end{description}
\ealgo

The use of the secant certainly mitigates the difficulty of computing the derivative $g'(y^k)$ accurately in the Newton method. 
Computing the secant by $(g(y^k)-g(y^{k-1}))/(y^k-y^{k-1})$ involves  the computation of the
 denominator $y^k-y^{k-1}$, which is 
 almost exact. To compute the numerator with high accuracy, 
 we still need to compute the function value $g(y^k)$ accurately at each iteration $(k=0,1,\ldots)$. 
The computations that are not carried out  with sufficient accuracy may affect the entire iteration. In particular, 
if the computed secant is much smaller than the real secant or the computed value of $g(y^k)$ is much larger than its real value at the $k$th iteration, then the next iterate $y^{k+1}$ could be smaller than $y^*$ and 
Algorithm~\ref{algorithm:secant} cannot be continued in a consistent manner. 

The stopping criteria in the APG method combined with 
Algorithm~\ref{algorithm:secant} should be carefully set up. 
We recall 
the discussion on the stopping criteria of the APG method combined with 
Algorithm~\ref{algorithm:bisection} in Section 2.3.  To approximate 
$(\widehat{\X}(y^k),\widehat{\Y}_1(y^k),\widehat{\Y}_2(y^k))$ 
with high accuracy %quality 
by the APG method, 
we modify the stopping criteria mentioned there.  
We introduce two criteria which must be satisfied simultaneously 
for the APG method to stop the iteration. 
The first one is $\kappa(\X,\Y_1,\Y_2) < tol$. 

For the second criterion, let 
$(\widetilde{\X}^p,\widetilde{\Y}_1^p,\widetilde{\Y}^p_2)$ denote the 
$p$th iterate of the APG method $(p=1,2,\ldots)$. We compute 
\begin{eqnarray*}
z^p = y^k - \alpha^k \frac{\|\widetilde{\X}^p\|(y^k-y^{k-1})}{\|\widetilde{\X}^p\| - \|\X^{k-1}\|} \ 
(p=1,2,\ldots,). 
\end{eqnarray*}
% Next, 
Here $z^p$ becomes the next iterate $y^{k+1}$ of Algorithm~\ref{algorithm:secant} if the APG method stops at  iteration $p$. 
Through numerical experiments, we have observed that the sequence $\{z^p\}$ increased 
on average as $p$ increase,  and in many cases, it  converged 
to some $z^*$.  %in many cases. 
To check the convergence of the generated sequence numerically  to some 
$z^*$, which is almost equal to $z^p$ at the $p$th iteration, 
some variance  of the subsequence $z^{q(p)+1},z^{p(q)+2},\ldots,z^{p}$ 
and the slope of their linear interpolation 
are incorporated into the second criteria (for example, 
$q(p) =  \max\{1, p-100\}$). Our  second criterion is: the variance and the slope 
are both sufficiently small. The details are omitted here. 

The second criterion can be very strict for some problems.
In fact, the APG method sometimes does not terminate with the two stopping criteria,
even when $p$ reaches a prescribed maximum number of iterations.
From our numerical experience, we noticed that this  usually occurs when the second criterion is 
too strict. 
For such cases, we allow Algorithm~\ref{algorithm:secant}  to switch  to the bisection temporarily. 
% when the APG method has not terminated within 
% the prescribed iteration upper bound. 
More precisely,  we add $\lb = \lbb^0$ and $\ub = y^0$ in Step 0 
and insert the following steps between Steps 2 and 3. \vspace{-2mm}
\begin{description}
\item{Step 2.3: } If $\|\X\| < \epsilon$, then let $\lb = \max\{\lb,y^k\}$. Otherwise let $\ub = \min \{\ub,y^k\}$.  
If the APG method for computing $(\X,\Y_1,\Y_2)$ at Step 1  
terminates  with 
the two stopping criteria, go to Step 3.  Otherwise go to Step 2.6. \vspace{-2mm}
\item{Step 2.6: } 
Reset $y^k = (\ub+\lb)/2$. If $\ub-\lb < \delta$, then 
stop the iteration, and output $\lbb^{k+1}$ as a valid lower bound for $y^*$ and $y^k$ as an 
approximation of $y^*$. Otherwise go to Step 1. % \vspace{-2mm}
\end{description}
This modification improves the robustness and stability of Algorithm~\ref{algorithm:secant}.

\section{Applications to DNN relaxations of Combinatorial Quadratic Optimization Problems}

Let $\C \in \SymMat^n$, $\cc \in \Real^n$, 
$\A \in \Real^{\ell \times n}$, $\b \in \Real^{\ell}$, 
% \red{$I_{\rm bin} \subset \{1,\ldots,r\}$}, 
$I_{\rm bin} \subset \{1,\ldots,n\}$ (the index set for binary variables),
and $I_{\rm comp} \subset \{(j,k): 1\leq j < k \leq n \}$ (the index set for pairs of complementary variables). 
%For simplicity of notation, we assume that $I_{\rm bin} = \{1,\ldots,q\}$ for some $q \geq 0$; if $q=0$ then 
%$I_{\rm bin} = \emptyset$. 
Consider a QOP of the following form:
\begin{eqnarray}
\zeta_{\mbox{\scriptsize QOP}} & = & \inf\left\{\u^T\C\u + 2\cc^T\u: 
\begin{array}{l}
\u \in \Real^n_+, \  \A\u - \b = \0, \\
 u_i(1-u_i) = 0 \ (i\in I_{\rm bin}), \\
 u_j u_k = 0 \ ((j,k) \in I_{\rm comp})
\end{array}
\right\}. 
\label{eq:QOP50}
\end{eqnarray}
We assume that the feasible region of QOP~\eqref{eq:QOP50} is nonempty.

To derive a doubly nonnegative (DNN) relaxation for QOP~\eqref{eq:QOP50} of the form  COP~\eqref{eq:COP1p}, 
%for QOP~\eqref{eq:QOP50}, 
we introduce some notation and symbols.
Let $\Real^{1+n}$ denote the $(1+n)$-dimensional Euclidean space 
of column vectors $\x=(x_0,x_1,\ldots,x_n)$, and 
$\Real^{1+n}_+$ $= $ 
$\left\{\x \in \Real^{1+n} : x_i\geq 0 \ \right.$ $\left. (i=0,\ldots,n)\right\}$ (the nonnegative orthant of $\Real^{1+n}$). Let $\SymMat^{1+n}$ denote the space of $(1+n) \times (1+n)$ 
symmetric matrices with row and column indices $i=0,1,\ldots,n$, 
$\SymMat^{1+n}_+$ the cone of positive semidefinite matrices in $\SymMat^{1+n}$, $\SymN^{1+n}$ the cone of nonnegative matrices 
in $\SymMat^{1+n}$ and $\DNN^{1+n} = \SymMat^{1+n}_+ \cap \SymN^{1+n}$ (the {\em DNN cone}). 
By definition, $\x\x^T \in \DNN^{1+n}$ for every 
$\x = (x_0,\u) \in \Real^{1+n}_+$. 

Let
\begin{eqnarray}
\left.
\begin{array}{lcl} 
\spaceV & = & \SymMat^{1+n}, \ \coneK_1 \ = \ \SymMat^{1+n}_+, \vspace{2mm} \\ 
\coneK_2 \ & = & \ \left\{\X \in \SymN^{1+n} : 
\begin{array}{ll}
X_{0i} = X_{i0} = X_{ii} \ (i\in I_{\rm bin}), \\
X_{jk} = X_{kj} = 0 \ ((j,k) \in I_{\rm comp})
\end{array}
\right\}, \vspace{2mm} \\ 
\coneK &=& \coneK_1 \cap \coneK_2 = \left\{\X \in \DNN^{1+n} : 
\begin{array}{ll}
X_{0i} = X_{i0} = X_{ii} \ (i\in I_{\rm bin}), \\
X_{jk} = X_{kj} = 0 \ ((j,k) \in I_{\rm comp})
\end{array}
\right\}, \vspace{2mm} \\ 
\H^0 &=& 
\mbox{the matrix in $\SymMat^{1+n}$ with $1$ at the $(0,0)$ element and $0$ elsewhere}, 
\vspace{2mm} \\
\Q^0 &= & 
\begin{pmatrix} 0 & \cc^T \\ \cc & \C \end{pmatrix}, \
\H^1 \ = \ \begin{pmatrix} -\b & \A \end{pmatrix}^T
\begin{pmatrix} -\b & \A \end{pmatrix} \in \SymMat^{1+n}_+ 
\subset \coneK_1^*+\coneK_2^*.
\end{array}
\right\} \label{eq:K1K2}
\end{eqnarray} 
Then, for every $\x = (x_0,\u) \in \Real^{1+n}_+$, 
\begin{eqnarray*}
%& & \x\x^T \in \DNN^{1+n}, \\ 
& & x_0 = 1 \ \Leftrightarrow \ \inprod{\H^0}{\x\x^T} = x^2_0 = 1, \\ 
& & \A\u-\b x_0 = \0 \ \Leftrightarrow \ \inprod{\H^1}{\x\x^T} = (\A\u-\b x_0)^T (\A\u-\b x_0) = 0, \\
& & \inprod{\Q^0}{\x\x^T} = \u^T\C\u + 2x_0\cc^T\u. 
\end{eqnarray*}
It is easy to verify that if $\x = (1,\u) \in \Real^{1+n}_+$, then 
\begin{eqnarray*} 
& & [\x\x^T]_{0i} = [\x\x^T]_{i0} = x_0u_i = u_i, \ [\x\x^T]_{ii} = u_i^2  \ (i \in I_{\rm bin} ), \\ 
& & [\x\x^T]_{jk} = 
u_ju_k = u_ku_j = [\x\x^T]_{kj}  \ ((j,k) \in I_{\rm comp});  
\end{eqnarray*}
hence,
\begin{eqnarray*}
 u_i(1-u_i) = 0 \ 
(i\in I_{\rm bin}) \ \mbox{ and } 
 u_j u_k = 0 \ ((j,k) \in I_{\rm comp})
%\end{array}
%\right\}
\ \Leftrightarrow \ \x\x^T \in \coneK_2
\end{eqnarray*}
Therefore, QOP~\eqref{eq:QOP50} is equivalent to 
\begin{eqnarray}
\zeta_{\mbox{\scriptsize QOP}} & = & \inf\left\{ \inprod{\Q^0}{\x\x^T} : 
\begin{array}{l}
\x\x^T \in \coneK_1 \cap \coneK_2, \\ [3pt]
\inprod{\H^0}{\x\x^T} = 1, \ \inprod{\H^1}{\x\x^T} =  0 
\end{array}
\right\}. \label{eq:QOP51}
\end{eqnarray}
Note that every feasible solution $\u \in \Real^n$ of QOP~\eqref{eq:QOP50}  
with the objective value $\u^T\C\u + 2\cc^T\u$ corresponds to a feasible solution 
$\x = (1,\u)$ of QOP~\eqref{eq:QOP51} with the same objective value $\inprod{\Q^0}{\x\x^T}$. By construction, 
it is obvious that $\H^0, \ \H^1 \in \SymMat^{1+n}_+ \subset \coneK_1^* + \coneK_2^*$; hence Condition (I) 
is satisfied.

Now, by replacing $\x\x^T \in \spaceV$ with a matrix variable $\X \in \spaceV$, 
we obtain COP~\eqref{eq:COP1p} with 
$\spaceV$, $\coneK_1$, $\coneK_2$, $\Q^0$, $\H^0$ and $\H^1$ given in~\eqref{eq:K1K2}, 
which serves as a DNN relaxation of QOP~\eqref{eq:QOP51}. 
For the efficient computation of the metric projection onto 
$\coneK_2$ defined above, 
we refer to the paper \cite{KIM2016}. 

\section{Preliminary numerical results}

We implemented the damped secant-bracketing method 
(Algorithm~\ref{algorithm:secant}) by modifying the software package BBCPOP \cite{ITO2017}, 
which is based on the BP method (Algorithm~\ref{algorithm:bisection}) for
solving DNN relaxations of QOPs and polynomial optimization problems (POPs) with binary, box and 
complementarity (BBC) constraints. For numerical tests, we experimented with  binary quadratic optimization 
problem (BQOP) instances from \cite{BIQMAC} in Section 6.1, 
and quadratic assignment problem (QAP) instances from \cite{QAPLIB} in Section 6.2. 
All the computations were performed in MATLAB 2018b on iMac Pro with Intel Xeon W CPU (3.2 GHZ) and 128 GB memory.

We compare our results from Algorithm~\ref{algorithm:secant} (the damped secant-bracketing method) with BBCPOP and SDPNAL+.
BBCPOP is based on Algorithm~\ref{algorithm:bisection} and SDPNAL+ \cite{YST2015} is a Matlab implementation of the majorized semismooth Newton-CG augmented Lagrangian method for large scale SDPs with bounded variables.

\subsection{Binary quadratic optimization problems}

We solved DNN relaxations of BQOP instances,  
bqp100-1,$\ldots$,bqp100-5,bqp500-1,$\ldots$,bqp500-5, from \cite{BIQMAC} by 
three methods, BBCPOP, 
% the secant-bracketing method based on 
Algorithm~\ref{algorithm:secant} and SDPNAL+ \cite{YST2015}.

Each BQOP instance is of the following form: 
\begin{eqnarray}
\zeta_{\rm BQOP} = \inf\left\{\v^T\F\v : \v \in \Real^r_+, \ v_i(1-v_i) = 0 \ (i=1,\ldots,r) \right\} \label{eq:BQOP50}
\end{eqnarray}
with $r=100$ or $r=500$, where $\F \in \SymMat^r$. 
The problem \eqref{eq:BQOP50} is a special case of QOP~\eqref{eq:QOP50}, thus %so that
the discussion in the previous section could be  directly applied to 
BQOP~\eqref{eq:BQOP50} for its DNN relaxation. To strengthen 
the relaxation \cite{KIM2013,KIM2017b}, instead of the direct application, we first reformulate BQOP~\eqref{eq:BQOP50} to the 
following QOP by introducing a slack variable vector $\w \in \Real^r_+$. 
\begin{eqnarray}
\zeta_{\rm BQOP} = \inf\left\{\v^T\F\v : 
\begin{array}{l}
 \v \in \Real^r_+, \ \w \in \Real^n_+, \ \v+\w = \e, \\
 v_i(1-v_i)=0, \ w_i(1-w_i)=1 \ (i=1,\ldots,r)  
\end{array} \right\}, \label{eq:BQOP51}
\end{eqnarray}
where $\e$ denotes the column vector of ones in $\Real^r$. 
Let $n = 2r$, $\ell = r$, $I_{\rm bin} = $\{1,\ldots,n\}, $I_{\rm comp} = \emptyset$, 
\begin{eqnarray*}
\u &=& \begin{pmatrix} \v \\ \w \end{pmatrix} \in \Real^n, \
\C  = \ \begin{pmatrix} \F & \O \\ \O & \O \end{pmatrix} \in \SymMat^n, \ \cc = \0 \in \Real^n, \ 
\A = \begin{pmatrix} \I & \I \end{pmatrix}, \ \b = \e, 
\end{eqnarray*}
where $\I$ denotes the $r \times r$ identy matrix. Then %we can rewrite
 BQOP~\eqref{eq:BQOP51} can be rewritten as QOP~\eqref{eq:QOP50}, and %derive 
its DNN relaxation of the form COP~\eqref{eq:COP1p} can be derived with 
$\coneK_1$, $\coneK_2$, $\Q^0$, $\H^0$ and $\H^1$ given in~\eqref{eq:K1K2}. 

We applied BBCPOP and Algorithm~\ref{algorithm:secant} with the modification mentioned at the end of Section 4 
to the Lagrangian relaxation of the resulting COP~\eqref{eq:COP1p} and its dual~\eqref{eq:COP1d}, {\it i.e.}, 
the pair of COPs~\eqref{eq:COP1Lp}  and~\eqref{eq:COP1Ld}  with $\lambda=10,000$. (See Section 2.2). 
Since $\coneK_1$ and $\coneK_2$ are given 
in~\eqref{eq:K1K2} with $I_{\rm bin} = \{1,\ldots,n\}$ and  $I_{\rm comp} = \emptyset$,  COP~\eqref{eq:COP1p} is equivalent to 
the DNN problem 
\begin{eqnarray}
\eta^p & = & 
\inf \; \left\{ \inprod{\Q^0}{\X} : 
\begin{array}{l}
\X \in \DNN^{1+n}, \ \inprod{\H^0}{\X} = 1, \ \inprod{\H^1}{\X} = 0, \\
X_{0i} = X_{i0} = X_{ii} \ (i=1,\ldots,n)  
\end{array}
\right\}. \label{eq:COP6p} 
\end{eqnarray}
We also applied SDPNAL+ to this problem with the default parameters. % setting. 

Table~\ref{table:BQOP} presents the numerical results. 
We can observe two distinct features of 
Algorithm~\ref{algorithm:secant}, which are inherited from the Newton-bracketing method to some extent, in comparison to BBCPOP. 
The first one is that `it' in Algorithm~\ref{algorithm:secant} is 
smaller than `it' in BBCPOP in most of the instances. This is 
due to the advantage of the secant iteration over the bisection 
iteration.  We note that Algorithm~\ref{algorithm:secant} sometimes switched to 
the bisection (Steps 2.3 and 2.6 stated at the end of Section 4)  in many instances, which increased
`it',  but `it' of  Algorithm~\ref{algorithm:secant} still was smaller than 
`it' of BBCPOP. 
However, `itAPG' (the total number of iterations spent by the APG method)  in Algorithm~\ref{algorithm:secant} is 
larger than `itAPG' 
in BBCPOP. 
This is because Algorithm~\ref{algorithm:secant} requires more accuracy 
in the 
computation of $(\widehat{\X}(y^k), \widehat{\Y}_1(y^k),\widehat{\Y}_2(y^k))$
 by the APG method as we 
discussed at the end of Section~4. As a result, the execution 
time of Algorithm~\ref{algorithm:secant} is longer than that of 
BBCPOP. 

The other important feature of Algorithm~\ref{algorithm:secant} is the quality of lower bounds obtained. 
Indeed, the lower bounds 
obtained by Algorithm~\ref{algorithm:secant} for larger scale instances are 
obviously tighter than those obtained by the other methods. 

%\FloatBarrier

\begin{table}[htp]
%\scriptsize{
\begin{center}
\caption{Binary QOP instances from \cite{BIQMAC}. 
LB denotes a valid lower bound for the optimal value. 
LB corresponds to $\lbb$ in 
Algorithms~\ref{algorithm:secant} and~\ref{algorithm:bisection} (BBCPOP). 
We also applied Lemma~\ref{lemma:LB1} to the output of SDPNAL+ 
and computed a valid lower bound for the optimal value. `sec' denotes the execution time,
`itAPG'  the total number of iterations spent in the APG method in the case of 
Algorithms~\ref{algorithm:secant} and~\ref{algorithm:bisection} and the total number 
of iterations in the case of SDPNAL+, and `it'  the total number of iterations in 
Algorithms~\ref{algorithm:bisection} and~\ref{algorithm:secant}. 
The bold digits in the 
column LB show the cases where Algorithm~\ref{algorithm:secant} computed tighter lower bounds than the other two methods.
} 
\label{table:BQOP}
\vspace{2mm}
% \scalebox{0.78}
\begin{tabular}{|l|c|l|l|l|}% {|l|l|l|l|l|l|l|l|l|l|}
\hline
%                                & Known               &          New lower bounds        & Old best lower bounds \\ % &                        & +SDPT3 \\
%Problem                   & upper bounds    &          LBv(timeBP,itAPGR/itBP) & by BBCPOP (Mittelmann) \\%& SDPLBv(sec) \\
%\hline
Problem                   & Opt.Val.    & Solver        & LB(sec,itAPG/it) \\ % & SDPLB(sec) \\
\hline
                 bqp100-1 &        -7970 &      Algorithm~\ref{algorithm:secant} & {\bf-8036}(1.1e2,14980/14) \\ % & -8036(0.0e0) \\ 
                 bqp100-1 &                  &            BBCPOP & -8039(5.7e1, 8453/20) \\% & \\ %-8039(0.0e0) \\ 
%                 bqp100-1 &                 &     SparsePOP & -8380(2.0e2, 2242) &                          \\ 
                bqp100-1 &              &    SDPNAL+ & -8050(5.4e1, 5741) \\%& \\ 
 \hline
                 bqp100-2 &       -11036 &      Algorithm~\ref{algorithm:secant} & {\bf-11036}(1.1e2,15310/14) \\ %& -11036(0.0e0) \\ 
                 bqp100-2 &                  &            BBCPOP & -11043(3.3e1, 5265/20) \\%& \\%-11043(0.0e0) \\ 
%                  bqp100-2 &                 &     SparsePOP & -11489(3.4e2, 2542) &                          \\ 
               bqp100-2 &              &    SDPNAL+ & -11039(2.9e2,23868)\\% & \\ 
 \hline
                 bqp100-3 &       -12723 &      Algorithm~\ref{algorithm:secant} & -12723(5.4e1, 8455/11) \\%& -12723(0.0e0) \\ 
                 bqp100-3 &                   &            BBCPOP & -12723(2.0e1, 3874/21) \\% & \\%-12723(0.0e0) \\ 
%                 bqp100-3 &                 &     SparsePOP & -13153(2.1e2, 2542) &                          \\ 
                 bqp100-3 &              &    SDPNAL+ & -12723(1.8e1, 3321) \\%& \\ 
\hline
                 bqp100-4 &       -10368 &      Algorithm~\ref{algorithm:secant} & -10368(6.0e1, 9105/11) \\%& -10368(0.0e0) \\ 
                 bqp100-4 &                  &            BBCPOP & -10368(2.2e1, 4124/14) \\%& \\%-10368(0.0e0) \\ 
%                 bqp100-4 &                 &     SparsePOP & -10731(2.0e2, 2500) \\%&                          \\ 
                 bqp100-4 &              &    SDPNAL+ & -10369(2.3e1, 3991) \\%& \\ 
\hline
                 bqp100-5 &        -9083 &      Algorithm~\ref{algorithm:secant} & -9083(1.0e2,13365/13) \\%& -9083(0.0e0) \\ 
                 bqp100-5 &                 &            BBCPOP & -9083(4.4e1, 7183/21) \\% & \\%-9083(0.0e0) \\ 
%                 bqp100-5 &                 &     SparsePOP & -9487(2.0e2, 2502) &                          \\ 
                 bqp100-5 &              &    SDPNAL+ & -9083(2.0e2,21039) \\%& \\ 
\hline
                 bqp500-1 &      -116586 &      Algorithm~\ref{algorithm:secant} & {\bf-122598}(1.3e4,35030/21) \\%& -122598(0.0e0) \\ 
                 bqp500-1 &                   &            BBCPOP & -133277(6.1e2, 2104/21) \\%& \\%-133277(0.0e0) \\ 
%                 bqp500-1 &                 &     SparsePOP & -125964(5.2e3, 2121) &                          \\ 
                 bqp500-1 &              &    SDPNAL+ & -122633(3.0e3, 8783) \\%& \\ 
\hline
                 bqp500-2 &      -128223 &      Algorithm~\ref{algorithm:secant} & {\bf-132730}(1.0e4,27840/17) \\% & -132730(0.0e0) \\ 
                 bqp500-2 &                   &            BBCPOP & -142163(5.6e2, 1915/21) \\%& \\% -142163(0.0e0) \\ 
%                bqp500-2 &                 &     SparsePOP & -136011(6.8e3, 2357) &                          \\ 
                bqp500-2 &              &    SDPNAL+ & -132773(3.1e3, 7591) \\%& \\ 
 \hline
                 bqp500-3 &      -130812 &      Algorithm~\ref{algorithm:secant} & {\bf-134800}(1.2e4,31350/18) \\%& -134800(0.0e0) \\ 
                 bqp500-3 &                   &            BBCPOP & -145715(7.0e2, 2323/22) \\% & \\%-145715(0.0e0) \\ 
%                 bqp500-3 &                 &     SparsePOP & -138453(5.3e3, 2120) &                          \\ 
                 bqp500-3 &              &    SDPNAL+ & -134841(3.1e3, 8656) \\%& \\ 
\hline
                 bqp500-4 &      -130097 &      Algorithm~\ref{algorithm:secant} & {\bf-135485}(1.0e4,28235/17) \\%& -135485(0.0e0) \\ 
                 bqp500-4 &                    &            BBCPOP & -147512(7.0e2, 2468/22) \\%& \\%-147512(0.0e0) \\ 
%                 bqp500-4 &                 &     SparsePOP & -139329(4.8e3, 2016) &                          \\ 
                 bqp500-4 &              &    SDPNAL+ & -135527(2.8e3, 7699) \\%& \\ 
\hline
                 bqp500-5 &      -125487 &      Algorithm~\ref{algorithm:secant} & {\bf-130300}(8.1e3,28630/18) \\%& -130300(0.0e0) \\ 
                 bqp500-5 &                   &            BBCPOP & -140946(7.1e2, 2355/22) \\% & \\%-140946(0.0e0) \\ 
%                 bqp500-5 &                 &     SparsePOP & -134092(1.2e4, 2422) &                          \\ 
                 bqp500-5 &                  &    SDPNAL+ & -130404(1.6e3, 6194) \\% & \\ 
\hline
\end{tabular}
\end{center}
%}
\end{table}

\subsection{Quadratic assignment problems}

We identify a matrix $\W = [\w_1,\dots,\w_r]\in \Real^{r\times r}$ 
with the $r^2$-dimensional column vector $\u = {\rm vec}(\W) =  [\w_1; \dots; \w_r] \in \Real^{r^2}$.
Given
matrices $\A, \B \in \Real^{r\times r}$, the quadratic assignment
problem is stated as 
\begin{eqnarray*}
\zeta_{\rm QAP} & = & \inf \left\{ 
\W \bullet \A\W\B^T : 
\begin{array}{l} 
 \e^T \W \e_j  =1, \ \e_j^T \W \e = 1, \\ [2pt]
 (j=1,\dots,r), \; \W \in \{0,1\}^{r\times r} 
 \end{array}
 \right\} \nonumber \\ [5pt]
& = &  \min \; \left\{ \u^T(\B\otimes \A) \u : 
\begin{array}{l}  
 (\e_j^T\otimes \e^T)\u  = 1, \ 
(\e^T\otimes \e_j^T)\u = 1\ \\ [3pt] 
 (j=1,\dots,r), \ % \u \in \{0,1\}^{r^2} 
u_i(1-u_i)=1 \ (i=1,\ldots,r^2)
\end{array} 
\right\}.
% \label{eq:QAP}
\end{eqnarray*}
Here $\e \in \Real^r$ denotes the vector of ones, $\e_j \in \Real^r$ the $j$th coordinate unit  vector, 
and $\otimes$ denotes the Kronecker product. This problem is a 
special cace of QOP~\eqref{eq:QOP50}, so  its 
DNN relaxation of the form COP~\eqref{eq:COP1p} can be derived. 

%\FloatBarrier

Table~\ref{table:QAP} shows numerical results on large scale QAP instances from \cite{QAPLIB}, 
whose optimal values are not known. COP~\eqref{eq:COP1p} derived as the 
DNN relaxation of the instance sko81 in Table~\ref{table:QAP} involves a 
$6562 \times 6562$ dense variable matrix $\X$. Such a large scale  DNN problem is extremely hard to solve 
by many of the existing software packages. To solve such a large scale QAP, it is 
meaningful to reduce the gap between the lower and upper bounds. 

 Recently,  Mittelmann \cite{MITTELMANN2018} applied BBCPOP to large scale 
 QAP instances with unknown optimal values, % including tai100a, tai100b and wil100, 
 and computed their then-new  lower bounds which  had not been 
 obtained
  before. 
His results demonstrated that BBCPOP is very powerful and effective for solving 
DNN problems from large scale QAP instances. His numerical experiments were conducted with % done 
MATLAB 2018a on 32GB   Intel(R) Core(TM) i7-7700K CPU @ 4.20GHz. To compare 
BBCPOP with Algorithm~\ref{algorithm:secant}, we applied the two methods to some of the QAP instances 
on the same computing environment (MATLAB 2018b on iMac Pro with Intel Xeon W CPU (3.2 GHZ) and 128 GB memory). See the 3rd and 4th columns of Table~\ref{table:QAP}.  
As observed in the numerical results on BQOP instances in Table~\ref{table:BQOP}, 
the lower bounds obtained by Algorithm~\ref{algorithm:secant} are of higher quality 
than those obtained by BBCPOP for all instances in Table~\ref{table:QAP},  but Algorithm~\ref{algorithm:secant} is $2$-$5$ times slower than BBCPOP. 
We also see that `it' of  Algorithm~\ref{algorithm:secant} is larger than `it' of BBCPOP 
for instances tai35b -- tai60b, tai80b and tho40. This is because the APG method often failed to terminate 
with the stopping criteria within the prescribed iteration upper bound, and switched to the 
bisection (Steps 2.3 and 2.6 stated at the end of Section 4) for those instances. 

In \cite{MITTELMANN2018}, BBCPOP was applied to larger scale QAP instances 
in addition to the ones in Table~\ref{table:QAP}, 
including tai100a, tai100b and wil100.  But we have not applied 
Algorithm~\ref{algorithm:secant} to those larger instances because they are too 
time consuming.

 \begin{table}[htp]
%\scriptsize{
\begin{center}
\caption{QAP instances from \cite{QAPLIB}.
LB denotes a valid lower bound for the optimal value. 
LB corresponds to $\lbb$ in Algorithms~\ref{algorithm:bisection} 
and~\ref{algorithm:secant} (BBCPOP). `sec' denotes the execution time, `itAPG' the total 
number of iterations spent in the APG method,  and `it'  the total number of iterations in  Algorithms~\ref{algorithm:secant} and BBCPOP. 
The bold  digits show differences between LBs obtained by 
Algortihm~\ref{algorithm:secant} and BBCPOP.
} 
\label{table:QAP}
% \scalebox{0.78}
\vspace{2mm}
\begin{tabular}{|l|l|l|l|l|l|}% {|l|l|l|l|l|l|l|l|l|l|}
\hline
                          & Best known       &         Algorithm~\ref{algorithm:secant}  & BBCPOP \cite{MITTELMANN2018}  \\ % &                        & +SDPT3 \\
Problem            & upper bounds    &          LB(sec,itAPG/it)    & LB(sec,itAPG/it)  \\%& SDPLBv(sec) \\
\hline
                   tai35b &    283315445 &       269{\bf 754270}(4.0e3,17265/28) & 269{\bf 532372}(1.3e3,6272/21) \\ 
                   tai40b &    637250948 &       60{\bf9242568}(9.6e3,19195/29) & 60{\bf8808404}(3.7e3,8168/21) \\ 
                   tai50b &    458821517 &       431{\bf245612}(2.4e4,15595/28) & 431{\bf090738}(7.3e3,5072/21)  \\ 
                   tai60a &      7205962 &       632{\bf6094}(2.8e4,6485/23)          & 632{\bf5979}(1.5e4,3710/22) \\ 
                   tai60b &    608215054 &      592{\bf829061}(6.6e4,16630/27) & 592{\bf371789}(3.5e4,9008/22) \\ 
%\hline
                  tai80a &      13499184    & 116570{\bf55}(1.1e5,5280/21)   & 116570{\bf14}(2.0e4,1091/22) \\ 
                  tai80b &      818415043 & 786{\bf900529}(3.7e5,18070/33) & {\bf298474}(8.1e4,4017/22) \\ 
%                   tai100a &     21044752 & Not solved & 17853840 \\ 
%                   tai100b &     1185996137 & Not solved & 1151591000 \\ 
% \hline
\hline
                    sko42 &        15812 &      1533{\bf4}(5.7e3,9245/16) & 1533{\bf3}(2.1e3,3898/21) \\ %1.533264e4 \\ 
                    sko49 &        23386 &       2265{\bf3}(1.2e4,9170/19) & 2265{\bf1}(4.4e3,3383/21) \\ % 2.265021e4 \\ 
                    sko56 &        34458 &      3338{\bf9}(2.3e4,8780/18) & 3338{\bf6}(8.7e3,3393/22) \\ % 3.338503e4 \\ %& 33389(2.8e2) \\ 
                    sko64 &        48498 &      470{\bf20}(5.6e4,8920/18)    & 470{\bf18}(1.8e4,3001/22) \\% 4.701738e4  \\                
                   sko72 &        66256 &  644{\bf61}(1.1e5,9565/17) & 644{\bf56}(5.3e4,4923/22) \\% 
%                sko72 &        66256 & Not solved & 64456(5.3e4,4923/22) \\% 6.445510e4  \\ 
                   sko81 &        90998 & 883{\bf70}(2.0e5,9555/17) & 883{\bf63}(8.1e4,3740/23) \\ 
%                   sko90 &       115534 & Not solved & 112424 \\
\hline
                    tho40 &       240516 &       226{\bf516}(5.0e3,10065/24) & 226{\bf491}(2.2e3,4826/21) \\ % 2.264901e5 \\ 
\hline
                    wil50 &        48816 &       4812{\bf3}(1.3e4,8650/20) & 4812{\bf2}(8.0e3,5453/21) \\ % 4.812102e4\\
%                    wil100 & 273038   &   Not solved & 268956 \\
\hline
\end{tabular}
\end{center}
\end{table}

\section{Concluding remarks}

While our focus has been on QOPs in this paper, 
the Lagrangian-DNN was extended to 
a class of polynomial optimization problems (POPs) with binary, box and complementarity 
constraints. In fact, BBCPOP was designed to solve such POPs, 
and the paper \cite{ITO2018} reported better numerical results in terms of computational time (and/or quality of lower 
bounds in some instances) on large scale 
randomly generated instances from the class than SDPNAL+ which is regarded as a state-of-art software for solving
DNN relaxation problems.

Theoretically, the extension of 
the Newton-bracketing method to the Lagrangian-DNN relaxation induced from the class 
of POPs is almost straightforward. 
As observed in Section 6, the quality of the optimal values of relaxation problems of large-scale QOP instances obtained by the secant-bracketing method
are clearly better than those by BBCPOP.
The secant-bracketing method is, however,  %much 
more time-consuming than 
BBCPOP.  As a result,  %it should be necessary
 the computational efficiency of the method should be improved before
%method first before 
incorporating it into BBCPOP to handle large scale QOPs and POPs.

%%%%%
%\input ref.tex

%\bibliographystyle{plain}
%\bibliography{./enhFOM}
%%%%%

\end{document}